%% file: ms.tex
\documentclass[letterpaper, 10 pt, conference]{ieeeconf}  

\IEEEoverridecommandlockouts                              

\usepackage{graphics} 
\usepackage{epsfig} 
\usepackage{amsmath} 
\usepackage{amssymb}  
\usepackage{booktabs}
\usepackage{units}
\usepackage{accents}
\usepackage{prettyref} \input{references.tex}
\usepackage{color}
\usepackage[super,negative]{nth}

\usepackage{tikz}
\usetikzlibrary{shapes,arrows}
\usepackage{pgfplots} 
\usepackage{pgfgantt}
\usepackage{pdflscape}
\pgfplotsset{compat=newest} 
\pgfplotsset{plot coordinates/math parser=false}
\newlength\fwidth
\usepackage{grffile}
\usetikzlibrary{plotmarks}
\usetikzlibrary{arrows.meta}
\usepgfplotslibrary{patchplots}
\usepackage{caption}
\usepackage{subcaption}
\usepackage{stmaryrd}
\usepackage{xcolor}
\usepackage{soul}

\usepackage{todonotes}

\newcommand{\mathmin}{\text{minimize}}
\newcommand{\mathst}{\text{subject to}}
\newcommand{\pplus}{\scriptscriptstyle +}

\newcommand{\T}{\scriptstyle\top}       

\newtheorem{assume}{Assumption}

\newtheorem{algo}{Algorithm}

\newtheorem{problem}{Problem}

\title{\LARGE \bf
A Stochastic Model Predictive Control Approach for Driver-Aided Intersection Crossing With Uncertain Driver Time Delay	
}

\author{Alexander Katriniok, Stefan Kojchev, Erjen Lefeber, Henk Nijmeijer
\thanks{A. Katriniok is with Ford Research \& Innovation Center, 52072 Aachen, Germany, {\tt\small de.alexander.katriniok@ieee.org}.}
\thanks{S. Kojchev, E. Lefeber and H. Nijmeijer are with the Dynamics and Control Group, Department of Mechanical Engineering, TU Eindhoven, 5600 MB Eindhoven, The Netherlands, {\tt\small S.Koychev@alumnus.tue.nl, A.A.J.Lefeber@tue.nl, H.Nijmeijer@tue.nl}. }}%

\begin{document}

\maketitle
\thispagestyle{empty}
\pagestyle{empty}

\begin{abstract}
We investigate the problem of coordinating human-driven vehicles in road intersections without any traffic lights or signs by issuing speed advices. The vehicles in the intersection are assumed to move along an \textit{a priori} known path and to be connected via vehicle-to-vehicle communication. The challenge arises with the uncertain driver reaction to a speed advice, especially in terms of the driver reaction time delay, as it might lead to unstable system dynamics. For this control problem, a distributed stochastic model predictive control concept is designed which accounts for driver uncertainties. By optimizing over scenarios, which are sequences of independent and identically distributed samples of the uncertainty over the prediction horizon, we can give probabilistic guarantees on constraint satisfaction. Simulation results demonstrate that the scenario-based approach is able to avoid collisions in spite of uncertainty while the non-stochastic baseline controller is not.
\end{abstract}

\input{content.tex}

\vspace*{-2mm}
\bibliographystyle{IEEEtran}
\bibliography{IEEEabrv}

\end{document}

%% file: references.tex
\newrefformat{sec}{\mbox{section \ref{#1}}}
\newrefformat{tab}{\mbox{Tab. \ref{#1}}}
\newrefformat{fig}{\mbox{Fig. \ref{#1}}}
\newrefformat{algo}{\mbox{Algorithm \ref{#1}}}

%% file: content.tex
\section{Introduction}
\label{sec:introduction}

Automating road intersections is a frequently discussed control problem, especially in the context of fully automated vehicles (AV) \cite{Chen2016}. The potential of automation to increase traffic flow, improve safety and reduce fuel consumption is significant. With a solution in place, one day we might even be able to dispense all traffic lights and signs. This contribution focuses on an intermediate solution, which aims at issuing speed advices to the driver, to achieve safe intersection crossing without any traffic lights or signs.

\subsection{Related Work}
\label{sec:introduction_relatedWork}

For the coordination of AVs in intersections, various solutions have been proposed, e.g., those based on multi-agent systems \cite{Kowshik2011}, hybrid system theory \cite{Hafner2013}, virtual platooning \cite{Medina2015} and model predictive control (MPC) \cite{Katriniok2017a,Campos2017,Shi2018}.

When discarding the assumption of fully automated vehicles, \cite{Schildbach2016} presents a robust MPC approach for determining safe gaps in the crossing traffic to pass the intersection and to optimize fuel efficiency. Thereby, no communication is available and only the human driven ego-vehicle can be controlled. Uncertainties in the motion of surrounding vehicles are covered by the robust approach. In \cite{Katriniok2017b}, the authors have proposed a distributed scenario-based MPC approach to orchestrate vehicles in intersections by issuing speed advices to the driver such that collisions between vehicles (or more generally agents) are avoided and traffic flow is optimized. For information exchange, the control scheme relies on vehicle-to-vehicle (V2V) communication. The driver is treated as an uncertain part of the control loop as his reaction to a speed advice might deviate from the expected one. This driver reaction is modeled as a proportional controller, which follows speed advices with a certain bounded offset, while the proportional gain and the speed offset are treated as uncertainties. However, the driver reaction time has been neglected in this first study.

\subsection{Main Contribution and Outline}
\label{sec:introduction_contribution} 

We extend our previous work in \cite{Katriniok2017b} by additionally accounting for an uncertain but constant driver reaction time, i.e., an uncertain but constant time delay. 
The control problem is to provide smooth driver speed advices for safe and efficient intersection crossing, even for an uncertain time delay. With the time delay, the open-loop prediction model might become unstable which is unfavorable in an MPC setting. In addition, the uncertain system response to a speed advice varies in a much wider range which complicates the calculation of smooth and convenient speed advices.

To solve the control problem, a scenario-based approach \cite{Schildbach2014,Farina2016} is pursued which draws independent and identically distributed (i.i.d.) samples from a bounded uncertainty set over the prediction horizon, referred to as scenarios. Essentially, every sample reflects a potential realization of the driver uncertainty. Ultimately, optimization is carried out over all scenarios subject to constraints that need to be satisfied for every scenario. With this methodology, we can give probabilistic guarantees on constraint satisfaction and eventually on collision avoidance. To account for unstable system dynamics, a state feedback gain is introduced which ensures stability for every uncertainty realization. Simulation results finally prove that the scenario-based approach is able to avoid collisions when the driver reaction is uncertain while the baseline MPC (neglecting uncertainty) is not. Hereafter, we mainly focus on a proof of concept while a real time solution of the control problem is part of ongoing research.

The paper is organized as follows. Section \ref{sec:modeling} outlines the MPC prediction model with a focus on the driver time delay extension. Section \ref{sec:OCP_feedbackGain} continues with the design of a feedback gain to prestabilize the MPC prediction model. Then, the distributed scenario-based MPC scheme is introduced in section \ref{sec:OCP} before \prettyref{sec:results} finally proves its efficacy. 

\subsection{Notation}
The predicted value of variable $x$ at the future time step $k+j$ is referred to as $x_{(k+j\mid k)}$. Moreover, $[x]_i$ refers to the $i$-th entry of vector $x$ while $\mathbb{N}^{\pplus}$ is the set of positive integers.

\section{Vehicle and Driver Reaction Model}
\label{sec:modeling}

To handle the complexity of intersection scenarios, we rely on the following assumptions:
\vspace*{1mm}
\begin{assume}[Intersection Scenarios]
A1. Only single intersection scenarios are considered; A2. The desired route of every agent is \textit{a priori} known; A3. All vehicles are human-driven; A4. Besides the driver reaction to a speed advice, there are no further uncertainties; \mbox{A5. All} vehicles are equipped with V2V communication; A6. No communication failures occur; \mbox{A7. Data} that has been transmitted after optimization at time $k$ is available to all other agents at time $k+1$. A8. Vehicle kinematic states are measurable. 
\end{assume}

\subsection{Vehicle Kinematics}
\label{sec:modeling_kinematics}

\begin{figure}[t!]	
	\begin{center}
	\def\svgwidth{7.0cm}	
	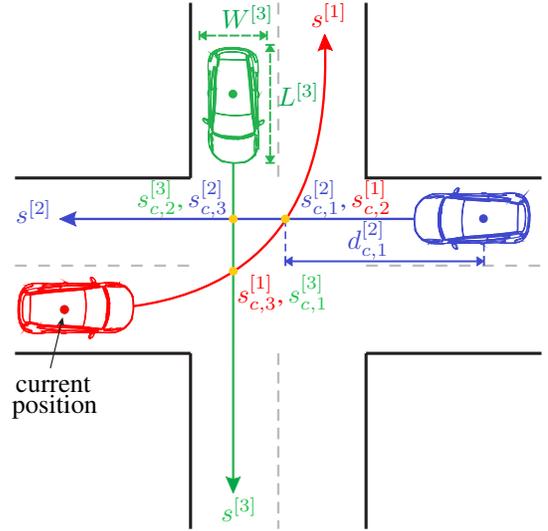
		\vspace*{-2mm}
	\caption{Schematic of the conflict resolution problem.}\vspace*{-8mm} 
	\label{fig:modeling_scheme}
	\end{center}
\end{figure}

Vehicle kinematics of every agent $i \in \mathcal{A}$ with \mbox{$\mathcal{A} \triangleq \{1, \ldots , N_A \}$} is formulated in terms of its acceleration $a_{x}^{[i]}$, velocity $v^{[i]}$ and path coordinate $s^{[i]}$ in the agent's reference frame with respect to the vehicle's geometric center, see \prettyref{fig:modeling_scheme}. The origin $s^{[i]}=0$ of agent $i$'s reference frame coincides with the first collision point $s_{c,l}^{[i]}$ with agent $l \in \mathcal{A}$ along his path coordinate $s^{[i]}$. In case, agent $i$ is not in conflict with any other agent, the origin refers to his initial position. The time evolution of velocity and position is represented as a double integrator while drivetrain dynamics are modeled as a first order lag element, i.e.,
\vspace*{-1mm}
\begin{align}
\frac{d}{dt}
\begin{bmatrix}
{a}_{x}^{[i]} \\
{v}^{[i]} \\
{s}^{[i]}
\end{bmatrix} &= 
\underbrace{\begin{bmatrix}
-\frac{1}{T_{a_x}^{[i]}} & 0 & 0 \\
                 1 & 0 & 0 \\
                 0 & 1 & 0
\end{bmatrix}}_{A_v^{[i]}}
\underbrace{
\vphantom{\begin{bmatrix}
\frac{1}{T_{a_x}^{[i]}} \\
0 \\
0
\end{bmatrix}}
\begin{bmatrix}
{a}_{x}^{[i]} \\
{v}^{[i]} \\
{s}^{[i]} 
\end{bmatrix}}_{x_v^{[i]}} +
\underbrace{\begin{bmatrix}
\frac{1}{T_{a_x}^{[i]}} \\
0 \\
0
\end{bmatrix}}_{B_v^{[i]}} \underbrace{\vphantom{\begin{bmatrix}
	\frac{1}{T_{a_x}^{[i]}} \\
	0 \\
	0
	\end{bmatrix}}
a_{x,\text{ref}}^{[i]}}_{u_v^{[i]}}
\label{eq:modeling_kinematics_ssModelVehicle}
\end{align}\\[-3mm]%
where $T_{a_{x}}^{[i]}$ denotes the dynamic drivetrain time constant and $a_{x,\text{ref}}^{[i]}$ the demanded acceleration.

\subsection{Driver Reaction Model}
\label{sec:modeling_driverModel}

We assume the driver to receive a speed advice $v_{\text{ref}}^{[i]}$ from the MPC controller and to translate this advice into a vehicle acceleration demand ${a}_{x,\text{ref}}^{[i]}$. Generally, we build upon our approach in \cite{Katriniok2017b}, in which the driver is modeled as a proportional controller with gain $K_d^{[i]} > 0$ and a bounded offset $\Delta v_{\text{d}}^{[i]}$ to the speed advice. The latter takes into consideration that the driver might not be able to exactly follow a speed advice. We extend this driver reaction model by means of a driver reaction time $\tau_d^{[i]} \geq 0$ to perceive the speed advice and react accordingly \cite{Treiber2006} --- corresponding to a time delay. The demanded vehicle acceleration ${a}_{x,\text{ref}}^{[i]}$, issued by the driver of agent $i$ as a reaction on the speed advice $v_{\text{ref}}^{[i]}$, can then be stated as 
\begin{align}
{a}_{x,\text{ref}}^{[i]}(t) = K_d^{[i]} \bigl[ v_{\text{ref}}^{[i]}(t-\tau_d^{[i]}) &+ \Delta v_{d}^{[i]}(t-\tau_d^{[i]}) \bigr.  \label{eq:modeling_driverModel_axref}\\ \bigl. &- v^{[i]}(t-\tau_d^{[i]}) \bigr].
\notag
\end{align}\\[-4mm]
During controller synthesis and system operation, we cannot be certain about the driver parameters. Therefore, we assume the driver model to be subject to an unmeasureable but bounded parametric uncertainty $\theta^{[i]} \triangleq (K_d^{[i]},\tau_d^{[i]})$ with
\vspace*{-1mm}
\begin{align}          
{K}_{d}^{[i]} \in [\underline{K}_{d}^{[i]},~ \overline{K}_{d}^{[i]}], ~~\tau_d^{[i]} \in [\underline{\tau}_d^{[i]},~  \overline{\tau}_d^{[i]}]. \label{eq:modeling_driverModel_uncertaintyParam}
\end{align}\\[-3mm]
Moreover, $\Delta v_{d}^{[i]}$ is treated as an unmeasurable but bounded additive uncertainty, i.e., 
\begin{align}          
\Delta v_{d}^{[i]} \in [\Delta\underline{v}_{d}^{[i]},~  \Delta\overline{v}_{d}^{[i]}] \label{eq:modeling_driverModel_uncertaintyAdd}
\end{align}
with $v_{\text{ref}}^{[i]} + \Delta v_{d}^{[i]}\geq 0$ to avoid negative speed advices.

\subsection{Resulting Prediction Model}
\label{sec:modeling_resultingModel}

The driver reaction-vehicle model with input and state delay 
is obtained when replacing $a_{x,\text{ref}}^{[i]}$ in \prettyref{eq:modeling_kinematics_ssModelVehicle} with \prettyref{eq:modeling_driverModel_axref}, i.e.,
\begin{align}
\dot{x}_v^{[i]}(t) = A_1^{[i]} {x}_v^{[i]}(t)  &+  A_2^{[i]} {x}_v^{[i]}(t-\tau_d^{[i]}) + B^{[i]} {v}_{\text{ref}}^{[i]}(t-\tau_d^{[i]}) \notag \label{eq:modeling_resultingModel_contBase} \\ &+ E ^{[i]}\Delta v_d^{[i]}(t-\tau_d^{[i]}) 
\end{align}
with $A_1^{[i]} \triangleq A_v^{[i]}$, $A_2^{[i]} \triangleq -K_d^{[i]} B_v^{[i]} [0~\, 1~\, 0]$, $B ^{[i]}\triangleq K_d^{[i]} B_v^{[i]}$ and $E^{[i]} \triangleq K_d^{[i]} B_v^{[i]}$. To be applied in the MPC framework, a discrete-time formulation of the model is required. In case of $\tau_d^{[i]}=0$, we discretize $(A_1^{[i]}+A_2^{[i]}, B^{[i]}, E^{[i]})$ using zero-order hold discretization as in \cite{Katriniok2017b}. 

For $\tau_d^{[i]}>0$, the continuous-time system \prettyref{eq:modeling_resultingModel_contBase} does not have a finite dimensional discrete-time representation \cite{Astrom1996}. We therefore \textit{digitalize} the driver, i.e., we assume the driver to sample the vehicle speed respectively the deviation from $v_{\mathrm{ref}}^{[i]}$ through a (digital) speedometer and to keep his acceleration demand constant between two sampling steps. This assumption translates in adding a zero-order hold element between the driver reaction \mbox{model \prettyref{eq:modeling_driverModel_axref}} and vehicle kinematics \prettyref{eq:modeling_kinematics_ssModelVehicle}, see \prettyref{fig:OCP_feedbackGain_controlLoop}. Thus, we are able to discretize both subsystems separately and gain a finite dimensional state space representation for any $\tau_d^{[i]} \in [\underline{\tau}_d^{[i]},\overline{\tau}_d^{[i]}]$. We define $\tau_d^{[i]} \triangleq \mathcal{T}^{[i]} \, T_s - \tilde{\tau}_d^{[i]}$ as an integer multiple $\mathcal{T}^{[i]} \in \mathbb{N}_0$ of the sampling time $T_s$ minus some remaining fraction \mbox{$0 \leq \tilde{\tau}_d^{[i]} < T_s$} of the time delay. With \prettyref{eq:modeling_driverModel_axref}, we can state the demanded acceleration $a_{x,\text{ref}}^{[i]}$ at the future time instance $t_{k+\mathcal{T}^{[i]}}$ in dependence of the delayed signals $v_{\text{ref}}^{[i]}$, $\Delta v_{\text{ref}}^{[i]}$ and $v^{[i]}$ at time $t_{k+\mathcal{T}^{[i]}} - \tau_d^{[i]} = t_k + \tilde{\tau}_d^{[i]}$, i.e.,
\begin{align}
a_{x,\text{ref}}^{[i]}(t_{k+\mathcal{T}^{[i]}}) \triangleq \, K_d^{[i]} \,&\Bigl[ v_{\text{ref}}^{[i]}(t_{k}+\tilde\tau_d^{[i]}) + \Delta v_{d}^{[i]}(t_{k}+\tilde\tau_d^{[i]}) \bigr. \notag\\ 
\Bigl.&-v^{[i]}(t_{k}+\tilde\tau_d^{[i]}) \Bigr]  \label{eq:modeling_resultingModel_axref}
\end{align}
\vspace*{-7mm}
\begin{align}
&\text{where~~~}v^{[i]}(t_{k}+\tilde\tau_d^{[i]}) \triangleq \Gamma_{A}^{[i]} \, x_v^{[i]}(t_{k}) + \Gamma_{B}^{[i]} \, a_{x,\text{ref}}^{[i]}(t_{k}), \notag\\
&\Gamma_{A}^{[i]} \triangleq
\begin{bmatrix}	0 &  1 & 0 \end{bmatrix} e^{A_{v}^{[i]}\tilde\tau_d^{[i]}}, ~ \Gamma_{B}^{[i]} \triangleq \begin{bmatrix}	0 & 1 & 0 \end{bmatrix} \int_{0}^{\tilde{\tau}_d^{[i]}} \hspace*{-2.5mm}e^{A_{v}^{[i]} s} ds\,B_{v}^{[i]}. \notag
\end{align}\\[-4mm]
Before $a_{x,\text{ref}}^{[i]}(t_{k+\mathcal{T}^{[i]}})$ is eventually applied to the vehicle, it is delayed by the driver reaction time for $\mathcal{T}^{[i]}-1$ time steps. By introducing the time delay states \mbox{$x_{\tau}^{[i]}(t_k) \triangleq [a_{x,\text{ref}}^{[i]}(t_{k+\mathcal{T}^{[i]}-1}),\ldots,a_{x,\text{ref}}^{[i]}(t_{k})]^{\T}$}, we can summarize those observations as 
\vspace*{-1mm}
\begin{align}
x_{\tau}^{[i]}(t_{k+1}) =  
\begin{bmatrix}
\prettyref{eq:modeling_resultingModel_axref} \\
\begin{bmatrix} I_{\mathcal{T}^{[i]}-1} & 0_{\mathcal{T}^{[i]}-1,1}\end{bmatrix} x_{\tau}^{[i]}(t_{k})
\end{bmatrix}
\end{align} \\[-3mm]
where $I_{\mathcal{T}^{[i]}-1}$ is the $(\mathcal{T}^{[i]}-1) \times (\mathcal{T}^{[i]}-1)$ identity matrix. 
With the discrete-time vehicle kinematics model
\mbox{\(
( \,\bar A_{v}^{[i]} \triangleq e^{A_{v}^{[i]} T_s}, ~ \bar B_{v}^{[i]} \triangleq \int_{0}^{T_s} e^{A_{v}^{[i]} s} ds\,B_{v}^{[i]} \,),
\)}
the time evolution of the vehicle states $x_{v}^{[i]}$ can be stated as
\vspace*{-1mm}
\begin{align}
x_{v}^{[i]}(t_{k+1}) = \bar{A}_{v}^{[i]} x_{v}^{[i]}(t_{k})  +  \bar{B}_{v}^{[i]} a_{x,\text{ref}}^{[i]}(t_{k}).
\end{align}\\[-4mm]
By augmenting the state vector with the time delay states \mbox{$x_{\tau,k}^{[i]} \triangleq x_{\tau}^{[i]}(t_k)$} and using the vehicle velocity as control output, the resulting discrete-time linear system $\Sigma_\theta^{[i]}$ evolves as
\begin{align}
{\Sigma}_{\theta}^{[i]} \triangleq \Bigg\{ \begin{array}{rl}
x_{k+1}^{[i]} \hspace*{-3mm}&= A_{\theta}^{[i]} x_{k}^{[i]} + B_{\theta}^{[i]} u_{k}^{[i]} + E_{\theta}^{[i]} w_{k}^{[i]} \vspace{1mm}\\
y_{k}^{[i]} \hspace*{-3mm}&= C^{[i]} x_{k}^{[i]}
\end{array} \label{eq:modeling_kinematics_sysModelDiscr}
\end{align}
with
\vspace*{-2mm}
\begin{align}
A_{\theta}^{[i]} &= \begin{bmatrix}
\bar A_{v}^{[i]} & ~\,0_{3, \mathcal{T}^{[i]}-1} & \bar B_{v}^{[i]} \\
-K_d^{[i]} \, \Gamma_A^{[i]}  & 0_{1, \mathcal{T}^{[i]}-1}   & -K_d^{[i]} \, \Gamma_B^{[i]} \\
0_{\mathcal{T}^{[i]}-1, 3}  &  I_{\mathcal{T}^{[i]}-1} & 0_{\mathcal{T}^{[i]}-1,1} 
\end{bmatrix},   \\
B_{\theta}^{[i]} &= \begin{bmatrix}
0_{3,1} \vspace*{2pt} \\ K_d^{[i]} \vspace*{2pt} \\ 0_{\mathcal{T}-1,1}
\end{bmatrix}, E_{\theta}^{[i]}=B_{\theta}^{[i]},~ C^{[i]} = [0~~1~~0~~0_{1, \mathcal{T}^{[i]}}] \notag
\end{align}
where $x_k^{[i]} \triangleq [x_{v,k}^{[i],{\T}},\,x_{\tau,k}^{[i],{\T}}]^{\T}$ refers to the state vector, \mbox{$u_k^{[i]} \triangleq v_{\text{ref}}^{[i]}(t_k+\tilde{\tau}_d^{[i]})$} to the control input, \mbox{$w_k^{[i]} \triangleq \Delta v_{d}^{[i]}(t_k+\tilde{\tau}_d)$} to the additive disturbance and \mbox{$y_k^{[i]}=v^{[i]}(t_k)$} to the system output. The attentive reader might have noticed that the initial condition $x_\tau^{[i]}(t_0)$ depends on \mbox{$v^{[i]}(t_0+\tilde\tau_d^{[i]}-n T_s)$} with $n=1,\ldots,\mathcal{T}^{[i]}$. Therefore, in the MPC implementation, we measure and store the velocity $v^{[i]}$ with a frequency which is high enough to obtain an appropriate initial condition. As the MPC is run with the fixed sample time $T_s$, we know that $v_{\text{ref}}^{[i]}(t_k+\tilde{\tau}_d^{[i]}-n T_s) = v_{\text{ref}}^{[i]}(t_k-n T_s)$ holds for $n \in \mathbb{N}_0$.

\subsection{Distances Between Agents}
\label{sec:modeling_vehDistances}

The distance between two agents $i,l \in \mathcal{A}$ is defined according to \cite{Katriniok2017a}. If two agents are potentially in conflict, they share a common collision point $s_{c,l}^{[i]}$ respectively $s_{c,i}^{[l]}$ along their respective path coordinate. Otherwise, we define $s_{c,l}^{[i]} = s_{c,i}^{[l]} = \infty$. This way, the distance between agent $i$ and $l$ is defined as the sum of the absolute distances $d_{c,l}^{[i]} \triangleq \lvert s^{[i]} - s_{c,l}^{[i]} \rvert$ and \mbox{$d_{c,i}^{[l]} \triangleq \lvert s^{[l]} - s_{c,i}^{[l]} \rvert$} to the agents' joint collision point, when existing, and infinite otherwise, i.e.,
\begin{align}
d_{l}^{[i]} = \begin{cases} \lvert s^{[i]} - s_{c,l}^{[i]} \rvert + \lvert s^{[l]} - s_{c,i}^{[l]} \rvert &  ,s_{c,l}^{[i]},s_{c,i}^{[l]} \neq \infty\\ \infty & ,\text{otherwise}. \end{cases} \label{eq:modeling_vehDistances_distanceDef}
\end{align}

\section{Stabilizing Feedback Gain}
\label{sec:OCP_feedbackGain}

\begin{figure}[t!]
	\centering
	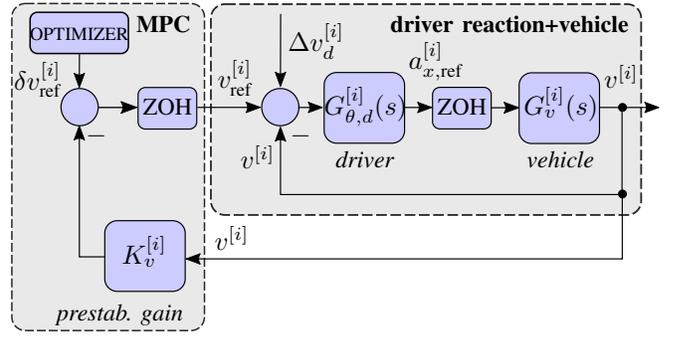
	\caption{Control loop with \textit{digitalized} driver, MPC controller and prestabilizing gain.}
	\vspace*{-7mm} 
	\label{fig:OCP_feedbackGain_controlLoop}
\end{figure}
Before further outlining the controller design, we briefly focus on how to handle unstable system dynamics, which might occur when a driver time delay is present. Without time delay, i.e., for $\tau_d^{[i]}=0$ it can easily be proven that the plant model with control input $v_{\mathrm{ref}}^{[i]}$ and output $v^{[i]}$ is strictly stable for every realization of the driver gain $K_d^{[i]} > 0$. When considering a driver reaction time $\tau_d^{[i]} > 0$, though, the prediction model might become unstable which is numerically unfavorable when it should be applied in the MPC framework \cite{Maciejowski2002}. In this case, a common approach is to design a prestabilizing state feedback gain $K_{\theta}^{[i]}$ and eventually apply the prestabilized plant model for prediction purposes \cite{Maciejowski2002}. This way, the control input can be written as 
\begin{align}
u_k^{[i]} = K_{\theta}^{[i]}  x_k^{[i]} + \delta u_k^{[i]}
\label{eq:OCP_feedbackGain_feedbackLaw}
\end{align} 
where $\delta {u}_k^{[i]} \triangleq \delta {v}_{\text{ref},k}^{[i]}$ is the new corrective control input that is determined by the MPC controller.

With the MPC as discrete-time controller, the feedback gain $K_{\theta}^{[i]}$ needs to be designed in the discrete-time domain as well. The main objective is to determine $K_{\theta}^{[i]}$ such that the closed-loop system $A_{\theta}^{[i]}+B_{\theta}^{[i]}K_{\theta}^{[i]}$ is Schur stable for all realizations of the uncertainty $\theta^{[i]}$. For this purpose, we implement a proportional feedback controller $K_v^{[i]}$ in accordance to \prettyref{fig:OCP_feedbackGain_controlLoop}. 

The discrete-time single-input single-output (SISO) \textit{driver reaction+vehicle} model in \prettyref{fig:OCP_feedbackGain_controlLoop} with control input $v_{\text{ref}}^{[i]}$ and control output $v^{[i]}$, can be represented in the z-domain as $G_\theta^{[i]}(z) = (G_{d,\theta}^{[i]}(z) \, G_{v}^{[i]}(z))/(1+G_{d,\theta}^{[i]}(z) \, G_{v}^{[i]}(z))$ where $G_{d,\theta}^{[i]}(z)$ and $G_{v}^{[i]}(z)$ denote the discrete-time counterparts of the continuous-time driver reaction and vehicle kinematics transfer functions.

\begin{figure}[h!]
	\centering	
	\setlength\fwidth{0.49\textwidth}
	\input{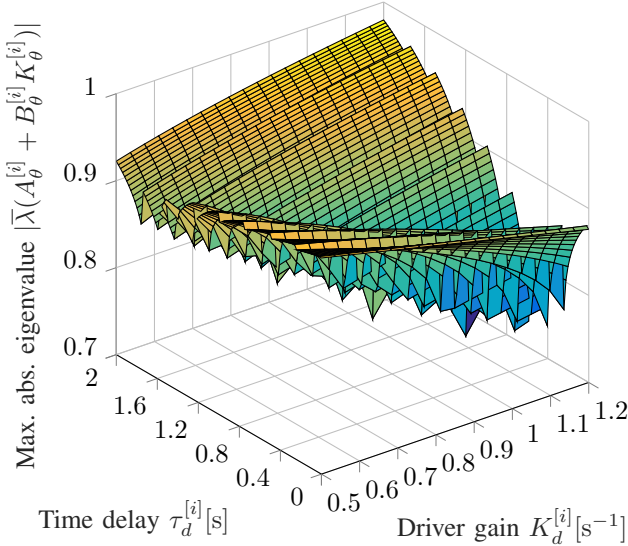}
	\setlength\fwidth{0.92\textwidth}\vspace*{-2mm}
	\caption{Max. abs. eigenvalue of the prestabilized system in \prettyref{sec:results} in dependence of $K_d^{[i]}$ and $\tau^{[i]}$ (2000 samples).}
	\vspace*{-2mm}
	\label{fig:OCP_feedbackGain_pzmap}
\end{figure}

We utilize the Nyquist criterion to design a proportional feedback gain $K_{v}^{[i]} < 0$ on the vehicle velocity to stabilize $G_\theta^{[i]}(z)$. Thus, we need to ensure that $1+K_{v}^{[i]} G_\theta^{[i]}(z)$ has only zeros inside the unit disc. It can be proven that if $G_\theta^{[i]}(z)$ is strictly stable, $G_\theta^{[i]}(z)/(1+K_{v}^{[i]} G_\theta^{[i]}(z))$ is strictly stable for all $K_{v}^{[i]} < 0$. If $G_\theta^{[i]}(z)$ is unstable, we obtain a lower bound $\underline{K}_{v}^{[i]}(\theta^{[i]})$ and an upper bound $\overline{K}_{v}^{[i]}(\theta^{[i]})$ on ${K}_{v}^{[i]}$ to ensure stability for the uncertainty realization $\theta^{[i]}$. 
When considering the parametric uncertainty  $\overline{\theta}^{[i]}=(\overline{K}_{d}^{[i]},\overline{\tau}_d^{[i]})$ with maximum gain and maximum time delay, we get the largest lower and smallest upper bound on ${K}_{v}^{[i]}$ for every possible realization of $\theta^{[i]}$. We choose 
$\underline{K}_{v}^{[i]}(\overline{\theta}^{[i]}) \leq {K}_{v}^{[i]} \leq \overline{K}_{v}^{[i]}(\overline{\theta}^{[i]})$
and finally obtain the state feedback gain $K_{\theta}^{[i]}$ as
\vspace*{-1mm}
\begin{align}
K_{\theta}^{[i]} \triangleq \bigl[ 0 ~~~ {-K_{v}^{[i]}} ~~~ 0 ~~~ 0_{1,\mathcal{T}^{[i]}} \bigr].
\label{eq:OCP_feedbackGain_Ktheta}
\end{align}\\[-4mm]
\prettyref{fig:OCP_feedbackGain_pzmap} illustrates for the system in \prettyref{sec:results} and 2000 samples of the uncertainty $\theta^{[i]}$ that the maximum absolute eigenvalue $\lvert\overline\lambda(A_\theta^{[i]}+B_\theta^{[i]} K_\theta^{[i]})\rvert$ of every prestabilized system is less than one, i.e., every prestabilized system is stable.
\vspace*{-2mm}

\section{Distributed Stochastic Optimal Control}
\label{sec:OCP}

\subsection{Control Problem and Algorithm}
\label{sec:OCP_problemAlgo}

The distributed control problem to solve can generally be stated as follows: 
\begin{problem}
In spite of the parametric and additive uncertainties \prettyref{eq:modeling_driverModel_uncertaintyParam}-\prettyref{eq:modeling_driverModel_uncertaintyAdd}, optimize the driver speed advice $v_{\text{ref}}^{[i]}$ for every agent $i$ such that collisions are avoided, the agent's cost function is minimized and constraints are satisfied.
\end{problem}
\vspace*{0.5mm}
The general idea of scenario-based MPC is to minimize, for every agent $i$, an average cost over scenarios $\kappa \in \mathcal{K}$, \mbox{$\mathcal{K} \triangleq \{1,\ldots,K\}$}, which are generated by drawing i.i.d. samples of the uncertainty over the prediction horizon, subject to constraints that need to be satisfied for every scenario \cite{Schildbach2014}. \prettyref{algo:OCP_problem_algoDist} sketches the algorithm which is applied to coordinate the agents in the intersection.
\begin{algo} \textbf{Scenario MPC at time $k$, Agent $i\in \mathcal{A}$}
	\begin{enumerate}
		\item \textbf{Receive data via V2V}: Receive distances  $d_{i,(\cdot\mid k)}^{[l]}$ to collision points from all agents $l \neq i$. 
		\item \textbf{Scenario Generation}: Sample $K$ scenarios.
		\item \textbf{Scenario Cost}: Determine an average scenario cost.
		\item \textbf{Scenario Constraints}: Impose input, state and safety constraints for every scenario.
		\item \textbf{Scenario Optimization}: Solve a single OCP which optimizes over $K$ scenarios s.t. scenario constraints. 
		\item \textbf{Broadcast data via V2V}: Broadcast distances $d_{l,(\cdot\mid k)}^{[i]}$.
		\item\textbf{Apply Control}:  Apply first element $u_{(k\mid k)}^{[i],\star}$. Go to 1).
	\end{enumerate}	\label{algo:OCP_problem_algoDist}
\end{algo}

\subsection{Scenario Model Predictive Control}
\label{sec:OCP_problem}
Hereafter, we outline the most important steps of \prettyref{algo:OCP_problem_algoDist}.

\subsubsection{\textbf{Scenario Generation}}
\label{sec:OCP_problem_scenarioGen}
During scenario generation, $K$ different scenarios are sampled, each representing a potential driver reaction in terms of the parametric uncertainty \mbox{$\theta^{[i]} = (K_d^{[i]},\tau_d^{[i]})$} and the additive uncertainty $\Delta v_d^{[i]}$. We assume that the driver does not change his general reaction over the prediction horizon. The deviation from the recommended speed, though, is considered to be time-varying over this interval, see \cite{Katriniok2017b}. To this end, we keep ${K}_{d}^{[i,\kappa]} \in [\underline{K}_{d}^{[i]},~ \overline{K}_{d}^{[i]}]$ and $\tau_d^{[i,\kappa]} \in [\underline{\tau}_d^{[i]},~  \overline{\tau}_d^{[i]}]$ constant over the prediction horizon for scenario $\kappa \in \mathcal{K}$ while the velocity offset $\Delta v_{d,(k+j\mid k)}^{[i,\kappa]}$ is sampled from the interval $[\Delta\underline{v}_d^{[i]},~  \Delta\overline{v}_d^{[i]}]$ for $\kappa \in \mathcal{K}$ and $j=0,\ldots,N-1$.

We eventually gain the following sampled system model for every scenario $\kappa \in \mathcal{K}$ (indicated by the superscript $\kappa$)
\begin{align}
{\Sigma}_{\theta}^{[i,\kappa]} \triangleq \Bigg\{ \begin{array}{rl}
x_{k+1}^{[i,\kappa]} \hspace*{-3mm}&= A_{\theta,K_\theta}^{[i,\kappa]} x_{k}^{[i,\kappa]} + B_{\theta}^{[i,\kappa]} \delta u_{k}^{[i]} + E_{\theta}^{[i,\kappa]} w_{k}^{[i,\kappa]} \vspace{1mm}\\
y_{k}^{[i,\kappa]} \hspace*{-3mm}&= C^{[i]} x_{k}^{[i,\kappa]}
\end{array} \raisetag{3.5mm}
\label{eq:OCP_problem_scenarioGen_sysModelSample}
\end{align}
where $A_{\theta,K_\theta}^{[i,\kappa]} \triangleq A_{\theta}^{[i,\kappa]}+B_{\theta}^{[i,\kappa]} K_\theta^{[i]}$ is the closed-loop matrix by applying feedback gain \prettyref{eq:OCP_feedbackGain_Ktheta}. Following \prettyref{eq:modeling_resultingModel_axref}, the time delay $\tau_d^{[i,\kappa]}$ is sampled from a continuous interval which eventually allows us to apply the stochastic MPC theory in \cite{Schildbach2014}.

Moreover, referring to \prettyref{eq:modeling_resultingModel_axref}, the initial condition of the time delay states $x_{\tau,0}^{[i]}$ depends on the uncertain speed offset $\Delta v_d^{[i]}$. This implies that we also need to sample the initial condition $x_0^{[i,\kappa]}$ as $x_{\tau,0}^{[i]}$ is neither measurable nor observable.

\subsubsection{\textbf{Scenario Cost}}
\label{sec:OCP_problem_scenarioCost}

The control objectives for every agent $i \in \mathcal{A}$ can be stated as follows: 
1) the velocity $v^{[i]}$ of every agent $i$ should follow the set speed $v_{\text{set}}^{[i]}$, being the minimum of the driver selected speed (usually the speed limit) and some situation dependent upper bound (e.g., in curves); 
2) driver speed advices should be smooth, as such step changes should be small; 
3) longitudinal accelerations should be minimized for efficient driving and 4) jerk should be minimized for the sake of comfort. We phrase these objectives in terms of the following quadratic objective function in dependence of scenario $\kappa \in \mathcal{K}$
\vspace*{-1mm}
\begin{align} 
J^{[i,\kappa]}(x_0^{[i,\kappa]},&\delta {u}_{(\cdot \mid k)}^{[i]}) \triangleq Q^{[i]} \sum^{N}_{j=1} (v_{\text{set},(k+j\mid k)}^{[i]} - v_{(k+j\mid k)}^{[i,\kappa]})^2 \notag \\
&+ R^{[i]} \sum^{N-1}_{j=0} \Delta {u}_{(k+j\mid k)}^{[i,\kappa],2}  \label{eq:nomOCP_local_costFcn} \\
&+ S_a^{[i]} \sum^{N}_{j=1} a_{x,(k+j\mid k)}^{[i,\kappa],2} + S_{\Delta a}^{[i]} \sum^{N}_{j=1} \Delta a_{x,(k+j\mid k)}^{[i,\kappa],2} \notag
\end{align} \\[-2mm]
where $x_0^{[i,\kappa]} = x_k^{[i,\kappa]}$ denotes the sampled initial condition  
at time $k$, \mbox{$\delta u_{(\cdot\mid k)}^{[i]} = [\delta u_{(k\mid k)}^{[i]}, \ldots, \delta u_{(k+N-1\mid k)}^{[i]} ]^{\T}$} the vector of corrective control actions over the prediction horizon of length $N$, \mbox{$\Delta {u}_{(k+j\mid k)}^{[i,\kappa]}={u}_{(k+j\mid k)}^{[i,\kappa]} - {u}_{(k+j-1\mid k)}^{[i,\kappa]}$} the step change of the resulting control input \mbox{$u_{(k+j\mid k)}^{[i,\kappa]} = K_{\theta}^{[i]}  x_{(k+j\mid k)}^{[i,\kappa]} + \delta u_{(k+j\mid k)}^{[i]}$} while $Q^{[i]} > 0$, $R^{[i]} > 0$, $S_a^{[i]} > 0$ and $S_{\Delta a}^{[i]} > 0$ are positive weights.

\subsubsection{\textbf{Scenario Constraints}} 
\label{sec:OCP_problem_cons}
Besides control objectives, local agent constraints and, most important, constraints for global collision avoidance need to be accommodated as well. 

In terms of \textbf{local agent constraints}, the speed advice $v_{\text{ref}}^{[i,\kappa]}$ plus speed offset $\Delta v_{d}^{[i,\kappa]}$ should be constrained for every agent $i \in \mathcal{A}$ such that only positive speeds are recommended and 
a driver selected upper speed bound (close to the speed limit) is accounted for. This claim is formulated as a constraint on the resulting control input $u_{(k+j\mid k)}^{[i,\kappa]}$ for \mbox{$j = 0,\ldots, N-1$} and every scenario $\kappa \in \mathcal{K}$, i.e.,
\begin{align} 
&u_{(k+j\mid k)}^{[i,\kappa]} \in \mathcal{U}_{(k+j\mid k)}^{[i,\kappa]} ~\forall \kappa \in \mathcal{K},~\text{with} \label{eq:OCP_problem_cons_inputCons} \\ 
&\mathcal{U}_{(k+j\mid k)}^{[i,\kappa]} \triangleq \Bigl\{ u \in \mathbb{R} \mid 0 \leq u + \Delta v_{d,(k+j\mid k)}^{[i,\kappa]} \leq \overline{v}_{(k+j\mid k)}^{[i]}  \Bigr. \notag \\
&\Bigl. ~~~~~~~~~~~~ 
\land \, u = K_\theta^{[i]} x_{(k+j\mid k)}^{[i,\kappa]} + \delta u_{(k+j\mid k)}^{[i]},~\delta u_{(k+j\mid k)}^{[i]} \in \mathbb{R} \Bigr\} \notag 
\end{align} \\[-4mm]
where $\overline{v}_{(k+j\mid k)}^{[i]}$ is an appropriately selected upper bound. Through the state constraint 
\begin{align} 
\hspace*{-1.5mm}
x_{(k+j\mid k)}^{[i,\kappa]} \in \mathcal{X}_{(k+j\mid k)}^{[i]} \triangleq \Bigl\{ x \in \mathbb{R}^3 & \mid~  \underline{a}_{x}^{[i]} \leq [x]_1 \leq \overline{a}_{x}^{[i]} 
\label{eq:OCP_problem_cons_stateCons} \\
& \land ~~\, 0 \leq [x]_2 \leq \overline{v}_{(k+j\mid k)}^{[i]} \Bigr\} \notag 
\end{align} \\[-4mm]
for $j=1,\ldots,N$ and every scenario \mbox{$\kappa \in \mathcal{K}$}, we also bound the actual velocity and constrain the vehicle acceleration to accommodate physical vehicle limitations and safe driving.

Following \cite{Katriniok2017b}, we impose the lower bound $\underline{v}_{\text{mean}}^{[i]}$ on the mean velocity over the prediction horizon when approaching a certain distance to the intersection, i.e., 
\begin{align} 
\frac{1}{N+1} \Bigl(v_{k}^{[i]} + \sum^{N}_{j=1} v_{(k+j\mid k)}^{[i,\kappa]} \Bigr) \geq \underline{v}_{\text{mean}}^{[i]},~\forall\kappa \in \mathcal{K}. 
\label{eq:OCP_problem_cons_minMeanVCons}
\end{align}  
Particularly, we claim that the prediction horizon at least covers the coordinate interval where potential collisions might occur with other agents. This way, convergence and feasibility of the distributed control scheme is ensured \cite{Katriniok2017b}. 

To ascertain \textbf{collision avoidance}, we follow our approach in \cite{Katriniok2017b}. Essentially, collision avoidance constraints need to be satisfied jointly, thus requiring a certain consensus among agents. To enforce consensus, we introduce time-invariant priorities on the agents that are determined once and held constant during the maneuver. Therefore, we define an injective prioritization function $\gamma: \mathcal{A} \rightarrow \mathbb{N}^+$ which assigns a unique priority to every agent, where a lower value corresponds to a higher priority. We specify the prioritized conflict set 
\(
\mathcal{A}_{c,\gamma}^{[i]} \triangleq \bigl\{  l \in \mathcal{A} \mid l \neq i \land \gamma(l) < \gamma(i) \land s_{c,l}^{[i]} \neq \infty \bigr\}
\) containing the agents $l \in \mathcal{A}$ which have a joint collision point with agent $i$ but a higher priority. Safety constraints can thus be phrased as
\begin{align}
d_{l,(k+j\mid k)}^{[i,\kappa]} \geq d_{\text{safe},l,(k+j\mid k)}^{[i]},~ \forall l \in \mathcal{A}_{c,\gamma}^{[i]}
\label{eq:OCP_problem_cons_safetyAbs}
\end{align}
for $j=1,\ldots,N$ and $\kappa \in \mathcal{K}$ with an appropriate safety distance $d_{\text{safe},l,(k+j\mid k)}^{[i]}$. Ultimately, only the agent with lower priority has to impose this safety constraint. With definition \prettyref{eq:modeling_vehDistances_distanceDef} of $d_{l,(k+j\mid k)}^{[i]}$, we can recast \prettyref{eq:OCP_problem_cons_safetyAbs} in the form \cite{Katriniok2017b}
\begin{align}
(s_{(k+j\mid k)}^{[i,\kappa]} - s_{c,l}^{[i]})^2 &\geq (d_{\text{safe},l,(k+j\mid k)}^{[i]} - d_{c,i,(k+j\mid k)}^{[l]})^2,~\label{eq:OCP_problem_cons_safetySquare} \\
&\forall l \in \mathcal{A}_{c,\gamma}^{[i]}: d_{\text{safe},l,(k+j\mid k)}^{[i]} > d_{c,i,(k+j\mid k)}^{[l]}. \notag
\end{align}
To avoid the necessity to transmit the trajectories of every scenario, 
every agent $l$ computes the distance  $d_{c,i,(k+j\mid k)}^{[l]}$ to its collision point based on the center of the interval $[\min_{\kappa \in \mathcal{K}} \{s_{(k+j\mid k)}^{[l,\kappa]}\}, \max_{\kappa \in \mathcal{K}} \{s_{(k+j\mid k)}^{[l,\kappa]}\}]$. The length of this interval, denoted as $\Delta L_{(k+j \mid k)}^{[l]}$, is leveraged to increase the safety distance of agent $i$. In the end, only $d_{c,i,(k+j\mid k)}^{[l]}$ and $\Delta L_{(k+j \mid k)}^{[l]}$ need to be transmitted to the other agents \cite{Katriniok2017b}.

\subsubsection{\textbf{Scenario Optimization}} 
\label{sec:OCP_problem_opt}

We have decomposed the control problem by separating the local cost functions and constraints while collision avoidance constraints are only imposed on agents with lower priority. Summarizing, the local OCPs can be cast as

\vspace*{2mm}
\textbf{Distributed Scenario OCP}, every agent $i \in \mathcal{A}$ solves:
\vspace*{-1mm}
\begin{align}
\underset{\delta {u}_{(\cdot \mid k)}^{[i]}}{\mathmin} ~& ~\, \frac{1}{K} \sum_{\kappa=1}^{K} J^{[i,\kappa]}(x_{0}^{[i,\kappa]}, \delta {u}_{(\cdot \mid k)}^{[i]}) \label{eq:OCP_problem_opt_OCPstatement}\\
\mathst~&~\, \text{system dynamics~} \text{(\ref{eq:OCP_problem_scenarioGen_sysModelSample})} \notag \\
&~\,\text{safety constraints~} \text{(\ref{eq:OCP_problem_cons_safetySquare})} \notag \\
&~\,\text{input (\ref{eq:OCP_problem_cons_inputCons}) \& state constraints~} \text{(\ref{eq:OCP_problem_cons_stateCons}),\,(\ref{eq:OCP_problem_cons_minMeanVCons})}. \notag 
\end{align}
Scenario OCP \prettyref{eq:OCP_problem_opt_OCPstatement} is a non-convex quadratically constrained quadratic program (QCQP). To solve the QCQP, we apply the penalty convex-concave procedure \cite{Lipp2016} as local method which iteratively solves a convex quadratic problem (QP). Finally, the control input $u_k^{[i],\star}=K_\theta^{[i]} x_k^{[i,1]} + \delta u_{(k\mid k)}^{[i],\star}$ is applied to the plant where $K_\theta^{[i]} x_k^{[i,1]} = \ldots = K_\theta^{[i]} x_k^{[i,K]}$.

\subsection{Constraint Violation Probability}

For a centralized scenario MPC scheme, \cite{Schildbach2014,Fagiano2015} have proven that scenario constraints are satisfied with a certain probability that depends on the number $K$ of scenarios. In this work, the uncertainties of every agent $i \in \mathcal{A}$ are assumed to be independent from each other, such that sampling can be carried out independently as well. Consequently, the theory in \cite{Schildbach2014} also holds for our distributed setup. According to \cite{Schildbach2014}, the most relevant criterion to guarantee closed-loop constraint satisfaction with a certain probability is the first predicted step constraint violation probability at time step $k+1$. With the parametric uncertainty $\theta^{[i]}$, the additive (uncertain) disturbance $\Delta v_d^{[i]}$ and a control input vector of dimension one, we obtain the following upper bound on the first predicted step (and as such closed-loop) state constraint violation probability for an arbitrary scenario $\kappa$
\vspace*{-1mm}
\begin{align}
\mathrm{P} \Bigl\{ x_{(k+1\mid k)}^{[i,\kappa]} \notin \mathcal{X}_{(k+1\mid k)}^{[i]} \Bigr\} \leq \frac{1}{1+K}.
\label{eq:OCP_consViol_PrViol}
\end{align} \\ [-4mm]
Given that only the worst case scenarios are broadcasted to the other agents (see \prettyref{sec:OCP_problem_cons}), collision avoidance constraints might be violated with an even lower probability. With the driver being eventually in charge of vehicle control, we consider a probabilistic guarantee on collision avoidance to be appropriate for the given application.

\section{Simulation Results}
\label{sec:results}
\begin{figure}[t!] 
	\begin{center}
		\includegraphics[width=4.4cm]{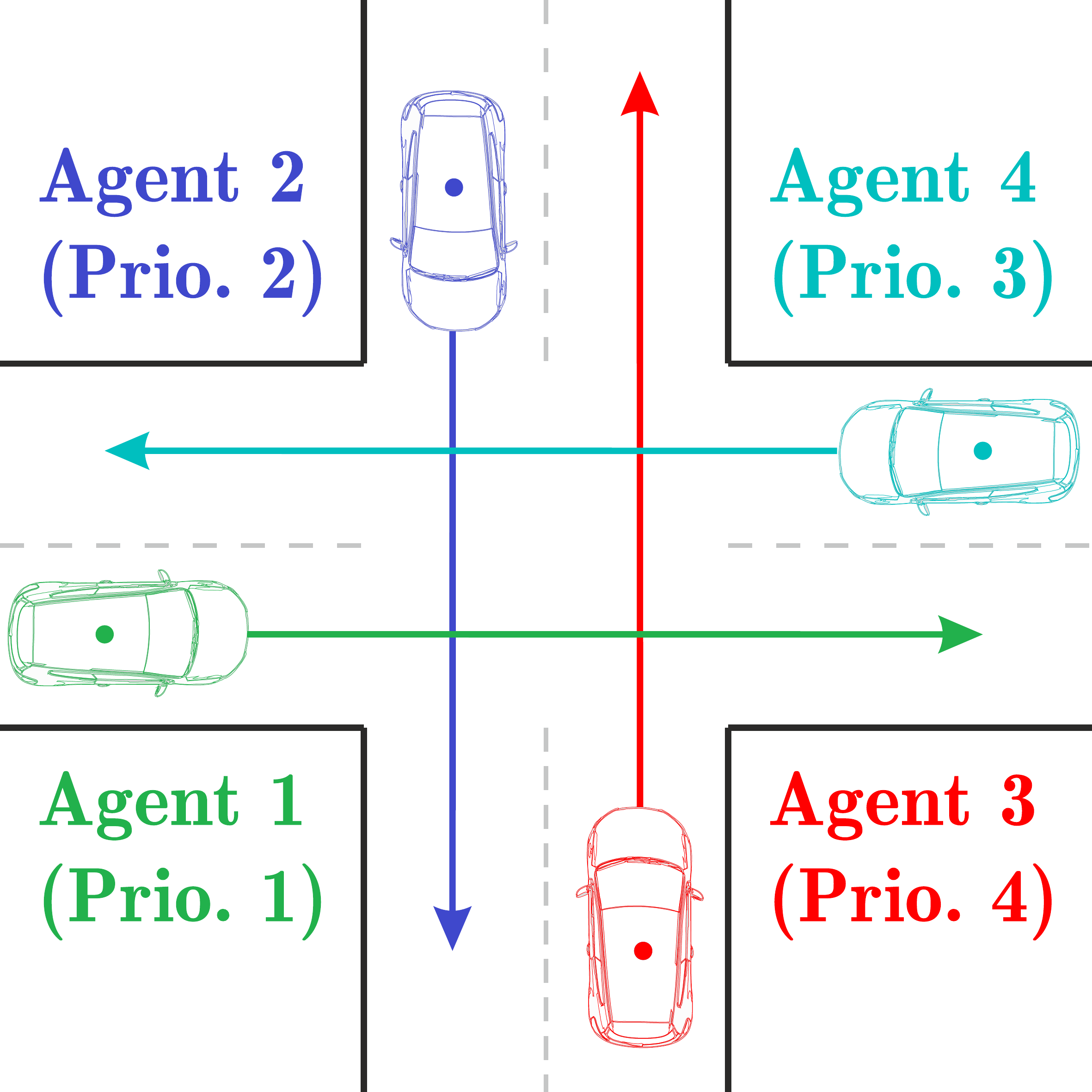}
				\vspace*{-2mm}
		\caption{Intersection scenario: Four straight passing agents.}\vspace*{-9mm}  
		\label{fig:results_setup_scenario}
	\end{center}
\end{figure}
\begin{figure*}[t!]
	\begin{center}
		\setlength\fwidth{0.82\textwidth}
		\input{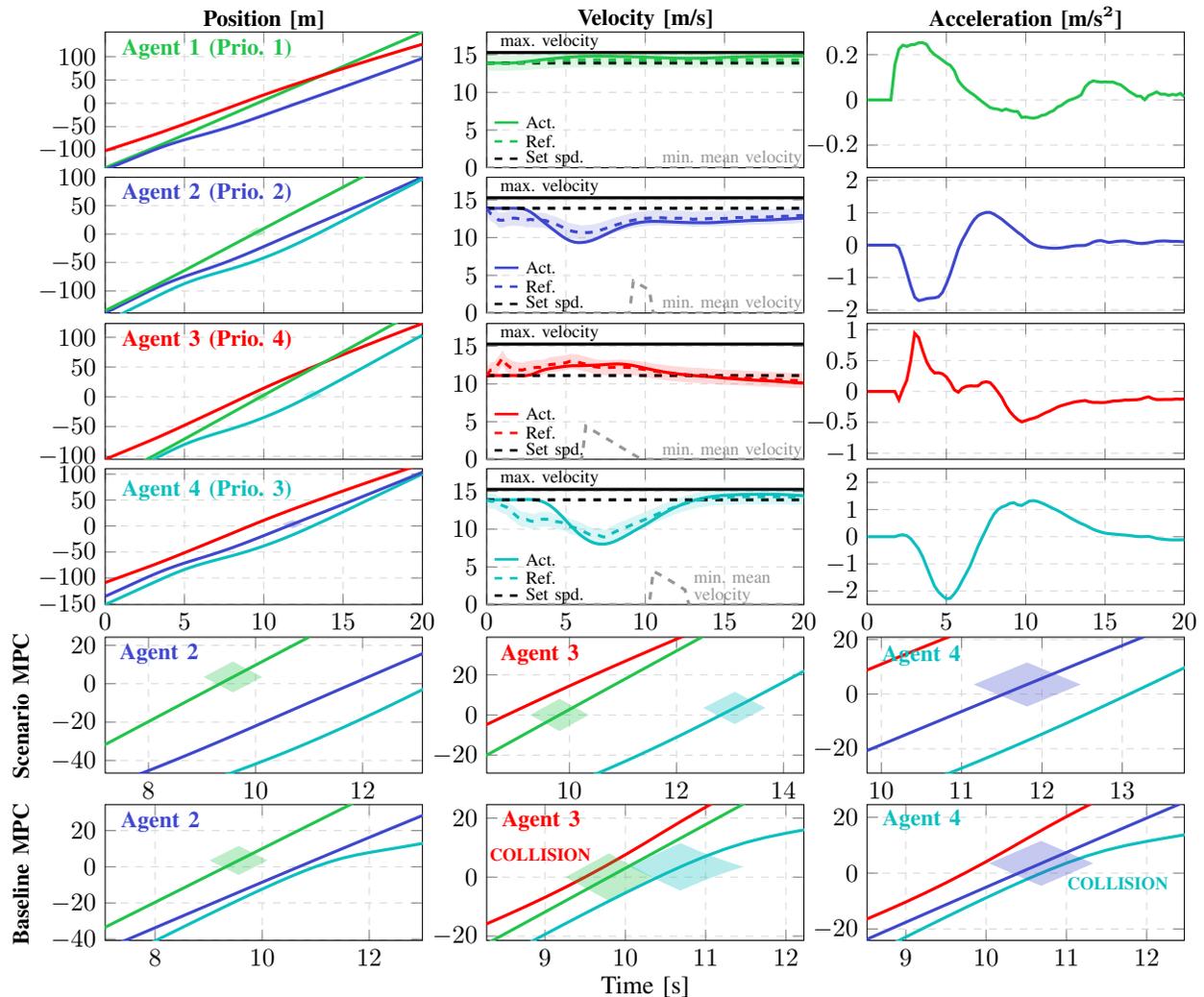}		
		\vspace*{-3mm}
		\caption{Scenario MPC is able to avoid collisions among agents in spite of uncertainty while Baseline MPC is not.} 
		\vspace*{-9mm}
		\label{fig:results_discussion_resultsPlot}
	\end{center}
\end{figure*}

\subsection{Simulation Setup}
\label{sec:results_setup}

To assess the validity of the stochastic approach, an urban four way intersection scenario with 
four agents passing the intersection straight is considered, see \prettyref{fig:results_setup_scenario}. Every agent has a length of $L^{[i]}=\unit[4.87]{m}$, a width of \mbox{$W^{[i]}=\unit[1.85]{m}$} and a dynamic drivetrain time constant of \mbox{$T_{a_x}^{[i]}=\unit[0.3]{s}$}. Parameter sets that represent a typical human driver are hard to determine. For the driver gain (in $\unit{s^{-1}}$) and the speed offset (in $\unit{m/s}$), we choose similar parameter intervals like in our previous study \cite{Katriniok2017b}. The time delay (in $\unit{s}$) is selected in accordance to studies on platooning \cite{Treiber2006}. In the simulation model, the following parameters are used: \mbox{$K_{d}^{[1]}=0.55$}, \mbox{$\tau_{d}^{[1]}=1.5$}, \mbox{$\Delta v_{d}^{[1]}=0.5$},
\mbox{$K_{d}^{[2]}=1.0$},  \mbox{$\tau_{d}^{[2]}=1.9$}, \mbox{$\Delta v_{d}^{[2]}=-0.3$}, \mbox{$K_{d}^{[3]}=0.7$},  \mbox{$\tau_{d}^{[3]}=1.8$}, \mbox{$\Delta v_{d}^{[3]}=-0.4$}, \mbox{$K_{d}^{[4]}=0.9$},  \mbox{$\tau_{d}^{[4]}=1.9$}, \mbox{$\Delta v_{d}^{[4]}=0.2$}. During simulation, $\Delta v_{d}^{[i]}$ is varied periodically in a step-wise manner by adding an offset which is bounded by the interval $[-0.2, 0.2]$.

For the stochastic MPC scheme, the uncertainty intervals are selected as: $\underline{K}_d^{[i]}=\unit[0.5]{s^{-1}}$, $\overline{K}_d^{[i]}=\unit[1.2]{s^{-1}}$ (driver gain bounds);   $\underline{\tau}_d^{[i]}=\unit[0]{s}$, $\overline{\tau}_d^{[i]}=\unit[2]{s}$ (time delay bounds);
$\Delta\underline{v}_d^{[i]}=\unit[-1]{m/s}$ and $\Delta\overline{v}_d^{[i]}=\unit[1]{m/s}$ (speed offset bounds).
Given those intervals, the feedback gain \mbox{$K_v^{[i]} \in [-1,\,-0.45]$} ensures a stable prediction model. As a trade-off between the resulting settling time and damping, we have chosen $K_v^{[i]}=-0.59$. Moreover, $K=99$ scenarios are generated for optimization --- with \prettyref{eq:OCP_consViol_PrViol}, this implies an upper constraint violation probability bound of $1 \%$.

The performance of the stochastic scheme is compared with a baseline MPC controller which only exploits a single realization of the driver parameters: 
\mbox{$K_{d}^{[1]}=0.6$}, \mbox{$\tau_{d}^{[1]}=1.2$}, \mbox{$K_{d}^{[2]}=0.8$},  \mbox{$\tau_{d}^{[2]}=1.5$}, \mbox{$K_{d}^{[3]}=0.9$},  \mbox{$\tau_{d}^{[3]}=1.2$}, \mbox{$K_{d}^{[4]}=1.1$},  \mbox{$\tau_{d}^{[4]}=1.0$} while $\Delta v_d^{[i]}$ is set to zero for every agent.

For both MPC regimes, the remaining parameters are set equally, i.e., using a sample time of \unit[0.25]{s} and a horizon length of \mbox{$N=40$} (i.e., a preview time of \unit[10]{s}). With this horizon length, the settling time is covered for $93 \%$ of all sampled system models (based on a 5000 sample analysis) which has turned out to be sufficient. Furthermore, the following weights are applied: \mbox{$Q^{[i]} = 0.5$}, \mbox{$R^{[i]} = 20$}, \mbox{$S_a^{[i]} = 5$}, \mbox{$S_{\Delta a}^{[i]} = 1$}. While agent 3 has selected a set speed of $\unit[11.1]{m/s}$, all other agents apply the speed limit of $\unit[13.9]{m/s}$ as set speed. The driver selected upper velocity bound $\overline{v}^{[i]}$ is set 10\% (i.e., $\unit[1.39]{m/s}$) higher than the speed limit. Longitudinal accelerations are bounded by \mbox{$\underline{a}_x^{[i]}= \unit[-7]{m/s^2}$} and $\overline{a}_x^{[i]}= \unit[4]{m/s^2}$. Moreover, agent priorities are fixed as in \cite{Katriniok2017b}: $\gamma(1)=1$, $\gamma(2)=2$, \mbox{$\gamma(3)=4$, $\gamma(4)=3$}. To solve OCP \prettyref{eq:OCP_problem_opt_OCPstatement}, qpOASES \cite{Ferreau2014} is utilized as QP solver.

\subsection{Discussion of Results}
\label{sec:results_discussion}

In \prettyref{fig:results_discussion_resultsPlot}, the $i$-th row illustrates the motion trajectories of \mbox{agent $i$} for the stochastic MPC scheme. The three respective columns, highlight from left to right: 1) the agent's path coordinate trajectory along with the trajectories of conflicting agents; \mbox{2) the} agent's actual (colored solid), maximum (solid black) and minimum mean velocity (dashed gray) together with the speed advice (colored dashed), the tolerated speed offset (colored patch) and the set speed (dashed black); 3) the actual vehicle acceleration. When agent $i$ imposes a safety constraint with respect to agent $l$, a colored polygon indicates the coordinate interval over time that must not be entered by the trajectory of agent $i$. 
The fifth row (scenario MPC) and sixth row (baseline MPC) show a closer insight into the time interval when agents are close to each other and collisions might occur. For reasons of brevity, we do not illustrate any motion trajectories for the baseline control scheme besides those in the sixth row. 

In the given scenario, agent 3, having the lowest priority, crosses the intersection first by speeding up from $\unit[11.1]{m/s}$ to $\unit[12.6]{m/s}$ with a moderate acceleration of $\unit[1]{m/s^2}$ to avoid collisions with agent 1 and agent 4. Evidently, the non-convex problem formulation bears the advantage to let agents cross in a sequence which is different from their priorities. After agent 3, agent 1 passes the intersection without the need to account for any other agent as he owns the highest priority. 
Agent 2 crosses the intersection after agent 1 by reducing his speed to $\unit[9.3]{m/s}$ with a maximum deceleration of $\unit[1.7]{m/s^2}$. Finally, \mbox{agent 4}, who exhibits the second lowest priority, crosses the intersection after agent 2. To safely avoid a collision with agent 2,
\mbox{agent 4} needs to slow down to $\unit[8]{m/s}$ with a maximum deceleration of $\unit[2.3]{m/s^2}$. For every agent, it is evident that the corresponding speed advices are very smooth and easy to follow for a human driver. Even despite a time delay of up to \unit[1.9]{s}, it can be recognized that the state and input trajectories satisfy constraints and do not show any noticeable oscillations. Without an appropriate feedback gain, simulation studies have shown that severe oscillations and unsmooth speed advices might occur. 

The last two rows in \prettyref{fig:results_discussion_resultsPlot} finally provide evidence that the scenario MPC scheme is able to avoid collisions between agents while the baseline MPC scheme, neglecting uncertainties, is not. Although agent 2 does not violate safety constraints for both strategies, agent 3 collides with agent 1 and agent 4 with agent 2 in case of the baseline controller. We can conclude that, despite uncertainties, the stochastic MPC scheme satisfies all our requirements.

\section{CONCLUSION AND FUTURE WORK}
\label{sec:conclusion}

We have proposed a distributed stochastic MPC scheme which provides speed advices to the driver in order to allow for safe and efficient intersection crossing without any traffic lights or signs. As an extension of our previous works, we include the driver reaction time delay as an additional parametric uncertainty in our control concept. Simulation results provide evidence that state, input and collision avoidance constraints are satisfied in spite of uncertainty. Future work aims at reducing the computational effort of the stochastic OCP and at verifying results in experimental tests.

%% file: intersection_3_w_turn_LW_Sc_new5.pdf_tex
\begingroup%
  \makeatletter%
  \providecommand\color[2][]{%
    \errmessage{(Inkscape) Color is used for the text in Inkscape, but the package 'color.sty' is not loaded}%
    \renewcommand\color[2][]{}%
  }%
  \providecommand\transparent[1]{%
    \errmessage{(Inkscape) Transparency is used (non-zero) for the text in Inkscape, but the package 'transparent.sty' is not loaded}%
    \renewcommand\transparent[1]{}%
  }%
  \providecommand\rotatebox[2]{#2}%
  \ifx\svgwidth\undefined%
    \setlength{\unitlength}{680.25bp}%
    \ifx\svgscale\undefined%
      \relax%
    \else%
      \setlength{\unitlength}{\unitlength * \real{\svgscale}}%
    \fi%
  \else%
    \setlength{\unitlength}{\svgwidth}%
  \fi%
  \global\let\svgwidth\undefined%
  \global\let\svgscale\undefined%
  \makeatother%
  \begin{picture}(1,1)%
    \put(0,0){\includegraphics[width=\unitlength,page=1]{intersection_3_w_turn_LW_Sc_new5.pdf}}%
	\definecolor{myRed}{rgb}{  1.00, 0.00, 0.00}%
	\definecolor{myGreen}{rgb}{0.1333333, 0.7694117, 0.2980392}%
	\definecolor{myBlue}{rgb}{ 0.2470588, 0.2823529, 0.8}%
	\definecolor{myCyan}{rgb}{ 0.05, 0.75, 0.75}%
	\put(0.29316178,0.6291634){\color[rgb]{0,0.65098039,0.31372549}\makebox(0,0)[lb]{\smash{}}}%
	\put(0.00169645,0.58099774){\color{myBlue}\makebox(0,0)[lb]{\smash{$s^{[2]}$}}}%
	\put(0.63169645,0.53199774){\color{myBlue}\makebox(0,0)[lb]{\smash{$d_{c,1}^{[2]}$}}}%
	\put(0.39138532,0.01491819){\color{myGreen}\makebox(0,0)[lb]{\smash{$s^{[3]}$}}}%
	\put(0.56215504,0.95414603){\color{myRed}\makebox(0,0)[lb]{\smash{$s^{[1]}$}}}%
	\put(0.39191238,0.95151305){\color[rgb]{0,0.65098039,0.31372549}\makebox(0,0)[lb]{\smash{$W^{[3]}$}}}%
	\put(0.49899409,0.80370992){\color[rgb]{0,0.65098039,0.31372549}\makebox(0,0)[lb]{\smash{$L^{[3]}$}}}%
	\put(0.00169645,0.26099774){\color[rgb]{0,0,0}\makebox(0,0)[lb]{\smash{current}}}%
	\put(-0.005169645,0.22099774){\color[rgb]{0,0,0}\makebox(0,0)[lb]{\smash{position}}}%
	
	\put(0.23169645,0.61099774){\makebox(0,0)[lb]{\smash{\color{myGreen}$s_{c,2}^{[3]}$,\,\color{myBlue}$s_{c,3}^{[2]}$}}}%
	\put(0.54169645,0.61099774){\makebox(0,0)[lb]{\smash{\color{myBlue}$s_{c,1}^{[2]}$,\,\color{myRed}$s_{c,2}^{[1]}$}}}%
	
	\put(0.42169645,0.43099774){\makebox(0,0)[lb]{\smash{\color{myRed}$s_{c,3}^{[1]}$,\,\color{myGreen}$s_{c,1}^{[3]}$}}}%
	
  \end{picture}%
\endgroup%

%% file: controlloop_ZOH_v3.pdf_tex
\begingroup%
  \makeatletter%
  \providecommand\color[2][]{%
    \errmessage{(Inkscape) Color is used for the text in Inkscape, but the package 'color.sty' is not loaded}%
    \renewcommand\color[2][]{}%
  }%
  \providecommand\transparent[1]{%
    \errmessage{(Inkscape) Transparency is used (non-zero) for the text in Inkscape, but the package 'transparent.sty' is not loaded}%
    \renewcommand\transparent[1]{}%
  }%
  \providecommand\rotatebox[2]{#2}%
  \ifx\svgwidth\undefined%
    \setlength{\unitlength}{338.76738703bp}%
    \ifx\svgscale\undefined%
      \relax%
    \else%
      \setlength{\unitlength}{\unitlength * \real{\svgscale}}%
    \fi%
  \else%
    \setlength{\unitlength}{\svgwidth}%
  \fi%
  \global\let\svgwidth\undefined%
  \global\let\svgscale\undefined%
  \makeatother%
  \begin{picture}(1,0.36573126)%
    \put(0,0){\includegraphics[width=\unitlength,page=1]{controlloop_ZOH_v3.pdf}}%
    \put(0.22972235,0.27084307){\color[rgb]{0,0,0}\makebox(0,0)[lb]{\smash{$v_{\text{ref}}^{[i]}$}}}%
    \put(0,0){\includegraphics[width=\unitlength,page=2]{controlloop_ZOH_v3.pdf}}%
    \put(0.12304028,0.07503584){\color[rgb]{0,0,0}\makebox(0,0)[lb]{\smash{$K_v^{[i]}$}}}%
    \put(0,0){\includegraphics[width=\unitlength,page=3]{controlloop_ZOH_v3.pdf}}%
    \put(0.31064532,0.21027675){\color[rgb]{0,0,0}\makebox(0,0)[lb]{\smash{$-$}}}%
    \put(0.08350886,0.20998343){\color[rgb]{0,0,0}\makebox(0,0)[lb]{\smash{$-$}}}%
    \put(0.65982539,0.27084257){\color[rgb]{0,0,0}\makebox(0,0)[lb]{\smash{$v^{[i]}$}}}%
    \put(0,0){\includegraphics[width=\unitlength,page=4]{controlloop_ZOH_v3.pdf}}%
    \put(0.0031611,0.2693605){\color[rgb]{0,0,0}\makebox(0,0)[lb]{\smash{$\delta v_{\text{ref}}^{[i]}$}}}%
    \put(0,0){\includegraphics[width=\unitlength,page=5]{controlloop_ZOH_v3.pdf}}%
    \put(0.0510626,0.01253126){\color[rgb]{0,0,0}\makebox(0,0)[lb]{\smash{$\small\textit{prestab. gain}$}}}%
    \put(0.42199179,0.33306284){\color[rgb]{0,0,0}\makebox(0,0)[lb]{\smash{$\small\textbf{driver reaction+vehicle}$}}}%
    \put(0,0){\includegraphics[width=\unitlength,page=6]{controlloop_ZOH_v3.pdf}}%
    \put(0.30806725,0.31185429){\color[rgb]{0,0,0}\makebox(0,0)[lb]{\smash{$\Delta v_{d}^{[i]}$}}}%
    \put(0.34829614,0.24105238){\color[rgb]{0,0,0}\makebox(0,0)[lb]{\smash{$G_{\theta,d}^{[i]}(s)$}}}%
    \put(0,0){\includegraphics[width=\unitlength,page=7]{controlloop_ZOH_v3.pdf}}%
    \put(0.14678741,0.24076411){\color[rgb]{0,0,0}\makebox(0,0)[lb]{\smash{$\small\text{ZOH}$}}}%
    \put(0.44247463,0.29180962){\color[rgb]{0,0,0}\makebox(0,0)[lb]{\smash{$a_{x,\text{ref}}^{[i]}$}}}%
    \put(0,0){\includegraphics[width=\unitlength,page=8]{controlloop_ZOH_v3.pdf}}%
    \put(0.47464016,0.24076411){\color[rgb]{0,0,0}\makebox(0,0)[lb]{\smash{$\small\text{ZOH}$}}}%
    \put(0,0){\includegraphics[width=\unitlength,page=9]{controlloop_ZOH_v3.pdf}}%
    \put(0.57033924,0.24105238){\color[rgb]{0,0,0}\makebox(0,0)[lb]{\smash{$G_{v}^{[i]}(s)$}}}%
    \put(0.13935053,0.33306284){\color[rgb]{0,0,0}\makebox(0,0)[lb]{\smash{$\small\textbf{MPC}$}}}%
    \put(0,0){\includegraphics[width=\unitlength,page=10]{controlloop_ZOH_v3.pdf}}%
    \put(0.22642023,0.09109056){\color[rgb]{0,0,0}\makebox(0,0)[lb]{\smash{$v^{[i]}$}}}%
    \put(0.25512854,0.18094995){\color[rgb]{0,0,0}\makebox(0,0)[lb]{\smash{$v^{[i]}$}}}%
    \put(0.36118794,0.18283891){\color[rgb]{0,0,0}\makebox(0,0)[lb]{\smash{$\small\textit{driver}$}}}%
    \put(0.57379299,0.18283877){\color[rgb]{0,0,0}\makebox(0,0)[lb]{\smash{$\small\textit{vehicle}$}}}%
    \put(0,0){\includegraphics[width=\unitlength,page=11]{controlloop_ZOH_v3.pdf}}%
    \put(0.02106389,0.32430377){\color[rgb]{0,0,0}\makebox(0,0)[lb]{\smash{$\scriptsize\text{OPTIMIZER}$}}}%
  \end{picture}%
\endgroup%